\def\marker{\>\hbox{${\vcenter{\vbox{
					\hrule height 0.4pt\hbox{\vrule width 0.4pt height 6pt
						\kern6pt\vrule width 0.4pt}\hrule height 0.4pt}}}$}\>}

\documentclass[12pt]{article}
\usepackage{indentfirst}
\usepackage{epsfig}
\usepackage{fullpage}
\usepackage{amsmath,amsfonts,amssymb,amsthm,graphics,amscd}
\usepackage{pgf,tikz}
\usepackage{graphicx}
\usepackage{subfigure}
\usepackage{amsmath,amsfonts,amssymb,amsthm,graphics,amscd,pst-plot,pstricks}
\usepackage{float}
\usepackage{graphicx}

\usepackage{xcolor}
\usepackage{caption}
\newtheorem {Theorem}  {Theorem}
\newtheorem {Lemma}{Lemma}[section]
\newtheorem {Corollary}[Theorem]{Corollary}

\theoremstyle{definition}
\newtheorem{Definition}{Definition}

\newtheorem{Conjecture}{Conjecture}
\numberwithin{Theorem}{section}
\newtheorem{remark}{Remark}

\newcommand{\D}{\Delta}

\newcommand{\phiv}{\varphi}

\newcommand{\CC}{\mathcal{C}}

\usepackage[
pdfauthor={},
pdftitle={From Nearly Triangle-Free to Nearly Bipartite in Edge-Chromatic Critical Graphs},
pdfstartview=XYZ,
bookmarks=true,
colorlinks=true,
linkcolor=blue,
urlcolor=blue,
citecolor=blue,
bookmarks=false,
linktocpage=true,
hyperindex=true
]{hyperref}

\begin{document}
\begin{center}
{\Large \bf From Nearly Triangle-Free to Nearly Bipartite in Edge-Chromatic Critical Graphs}

Xuli Qi*,\,\, Chunhui Ge

\vspace{3mm}
\baselineskip=0.25in
{\it Department of Mathematics, Hebei University of Science and Technology,}  \\
{\it Shijiazhuang 050018, P.R. China}\\
\vspace{3mm}
\end{center}	
\footnotetext{
E-mail address: xuliqi@hebust.edu.cn}

\begin{abstract}
Let $G$ be a simple $\Delta$-critical graph of order $n$,
and suppose that all triangles of $G$ contain a common vertex
$a$; equivalently, $H:=G-a$ is triangle-free. We prove that, under several  degree conditions, $H$ is bipartite. Consequently, $G$ is nearly bipartite and overfull.
First, for every fixed $\varepsilon,\xi>0$, there exists
$n_0:=n_0(\varepsilon,\xi)$ such that, whenever $n\ge n_0$,
the conditions $\Delta(G)\ge(\frac{2}{5}+\varepsilon)n$ and
$\delta(H)\ge(\frac{1}{4}+\xi)n$ imply that $H$ is bipartite. Second, suppose that there exists an integer $h$ satisfying
$1\le h<\Delta$ and $2\Delta-1<h(h+1)$. With $s:=\Delta-h$, the conditions $4s+\delta(H)\ge2n-1$ and $s+2\delta(H)\ge n$ then imply that $H$ is bipartite. Third, if $\sigma:=\Delta-d_G(a)$ satisfies  $4\sigma+\delta(H)\ge2n-5$ and $\sigma+2\delta(H)\ge n-1$, then $H$ is bipartite. In particular, $H$ is bipartite whenever $\delta(H)\ge\frac{2n}{7}$ and $\Delta\ge\frac{3n}{7}+\min\{d_G(a),\lceil\sqrt n\rceil\}$.
The proofs combine average-degree estimates, local
degree-deficiency bounds, multi-fans, Kierstead paths, and a
bipartite extension lemma.

{\small \it Keywords:}
edge coloring; edge-chromatic critical graph; nearly triangle-free graph; nearly bipartite graph; Overfull Conjecture; multi-fan; Kierstead path
\end{abstract}

\section{Introduction}
Unless otherwise stated, all graphs considered in this
paper are finite and simple. We generally follow the notation and terminology of Stiebitz et al.~\cite{SSTF2012}. Let $G=(V(G), E(G))$ be a simple graph, where $V(G)$ and $E(G)$ denote its vertex set and edge set, respectively. We use $\Delta(G)$ and $\delta(G)$ to denote the maximum degree and minimum degree of $G$, and simply write $\Delta := \Delta(G)$ and $\delta := \delta(G)$ throughout the paper. For $u \in V(G)$, let $d_G(u)$ and $N_G(u)$ denote the degree and the neighborhood of $u$, respectively. The average degree of $G$ is denoted by  $\overline{d}(G)$. For $V'\subseteq V(G)$, let $G-V'$ be the subgraph obtained from $G$ by deleting the vertex set $V'$ together with all its incident edges. In particular, for a vertex $v \in V(G)$, we write $G-v$ for $G-\{v\}$. For each edge $e \in E(G)$, let $G-e$ denote the graph obtained by deleting $e$ from $G$.   For any positive integer $k$, let $[k] :=\{1, \dots, k\}$ denote the set of the first $k$ positive integers. A (proper) $k$-edge-coloring of $G$ is a map $\varphi : E(G) \to [k]$ that assigns a color $\varphi(e) \in [k]$ to each edge $e \in E(G)$ such that no two adjacent edges share the same color. Let $\mathcal{C}^k(G)$ denote the set of all $k$-edge-colorings of $G$. The {\bf chromatic index} $\chi'(G)$ is the least integer $k$ for which $\mathcal{C}^k(G) \neq \emptyset$.

A graph $G$ is {\bf edge-chromatic critical} if $\chi'(H) < \chi'(G)$ for every proper subgraph $H$ of $G$. Furthermore, $G$ is called a {\bf $\Delta$-critical} graph if it is edge-chromatic critical and $\chi'(G)=\Delta(G)+1$. An edge $e$ is {\bf critical} if $\chi'(G-e)<\chi'(G)$.  Note that every $\Delta$-critical graph has minimum degree at least two. A subset $X\subseteq V(G)$ is an independent set if no two vertices in $X$ are adjacent, and $G$ is {\bf bipartite} if $V(G)$ can be partitioned into two independent sets. A graph $G$ is {\bf triangle-free} if it contains no cycle of length three.
	
In the 1960s, Vizing~\cite{Vizing1965-1} and, independently, Gupta~\cite{Gupta1967} proved that $\Delta(G) \le \chi'(G) \le \Delta(G)+1$, which naturally partitions graphs into two classes. Following Fiorini and Wilson~\cite{FW1977}, a graph $G$ is of {\bf Class 1} if $\chi'(G)=\Delta(G)$; otherwise, it is of {\bf Class 2}. In 1965, Vizing~\cite{Vizing1965-1} proved the following fundamental result, known as Vizing's Adjacency Lemma: if $G$ is a graph of Class 2 with maximum degree $\Delta(G)$ and $e=xy$ is a critical edge of $G$, then $x$ has at least $\Delta(G)-d_G(y)+1$ neighbors with maximum degree $\Delta(G)$ in $V(G)\setminus \{y\}$.
If a graph $G$ satisfies $|E(G)|>\Delta(G)\lfloor \frac{|V(G)|}{2} \rfloor$, any proper edge coloring of $G$ requires at least $\Delta(G)+1$ colors; such graphs are {\bf overfull}. A subgraph $H$ of $G$ is called {\bf $\Delta(G)$-overfull} if $|E(H)|>\Delta(G)\lfloor \frac{|V(H)|}{2} \rfloor$. Holyer~\cite{Holyer1981} proved that determining whether an arbitrary graph is of Class 1 is NP-complete. However, a  conjecture proposed by Chetwynd and Hilton would imply that the chromatic index of any graph $G$ satisfying $\Delta(G)>\frac{|V(G)|}{3}$ can be computed in polynomial time. The conjecture of Chetwynd and Hilton~\cite{CHETWYND1986,CHETWYND1989} is stated as follows.

\begin{Conjecture}[Overfull Conjecture]
Let $G$ be a graph with $\Delta(G)>\frac{|V(G)|}{3}$. Then $G$ is of Class $2$ if and only if $G$ contains a $\Delta(G)$-overfull subgraph.
\end{Conjecture}

Chetwynd and Hilton~\cite{CHETWYND1989} verified the Overfull Conjecture for graphs $G$ with $\Delta(G)\ge |V(G)|-3$.  Ge, Gutin, and Qi~\cite{GGQ2026} confirmed the conjecture for graphs with maximum degree 4. Shan~\cite{Shan2024-2} proved that for any $0<\varepsilon\le\frac{1}{14}$, there exists a positive integer $n_0$ such that if $G$ is a graph with order $|V(G)|\ge n_0$ and $\D(G)\ge(1-\varepsilon)|V(G)|$, then the conjecture holds. Recently, Chen, McDonald, and Shan~\cite{CMS2026} established that, for any $0<\varepsilon<1$, the conjecture holds for all sufficiently large graphs satisfying $\Delta(G) \ge \frac{(1+\varepsilon)|V(G)|}{2}$. Since every Class 2 graph with maximum degree $\Delta$ contains a $\Delta$-critical subgraph,
the Overfull Conjecture is equivalently stated as follows: every $\Delta$-critical graph $G$ with $\Delta(G)>\frac{|V(G)|}{3}$ is overfull. In 2022, Cao, Chen, Jing, and Shan~\cite{CCJS2022} proved  that the conjecture is true for $\Delta$-critical graphs $G$ of order $n$ provided that $\Delta(G)-\frac{7\delta(G)}{4}\ge \frac{3n-17}{4}$. Qi, Ge, and Feng~\cite{QGF2025} improved this condition to $\Delta(G)-\frac{5\delta(G)}{3}\ge \frac{2n-7}{3}$.

We call a graph $G$ {\bf nearly triangle-free} if there exists a vertex $a\in V(G)$ such that $G-a$ is triangle-free; equivalently, all triangles of $G$ contain a common vertex. Similarly, $G$ is called {\bf nearly bipartite} if deleting at most one vertex results in a bipartite graph.

In this paper, we consider a $\Delta$-critical graph
$G$ with a specified vertex $a$ such that $H:=G-a$ is triangle-free. We give degree conditions forcing $H$ to be bipartite. Consequently, $G$ is nearly bipartite and
overfull. We establish one asymptotic criterion and two finite-order criteria. The asymptotic proof uses an average-degree estimate to find a large bipartite subgraph of $H$, and then extends its bipartition to all of $H$. The finite-order proofs instead rely on local degree-deficiency bounds, multi-fans, and Kierstead paths.

\begin{remark}\label{rem:max-degree}
Since $G$ is $\Delta$-critical, every edge of $G$
is critical. Choose a vertex $x \in N_G(a)$. Applying Vizing's Adjacency Lemma to an
edge $ax$ incident with $a$, we can choose $b\in N_G(a)\setminus\{x\}$ such that $d_G(b)=\Delta$. Applying the lemma again to the critical edge $ab$, with $b$ as the first endvertex, we can choose $c\in N_G(b)\setminus\{a\}$ such that
$d_G(c)=\Delta$. Since $H=G-a$, we have $d_H(b)=\Delta-1$ and $d_H(c)\ge\Delta-1$. Moreover, since $H$ is triangle-free and
$bc\in E(H)$, we have $N_H(b)\cap N_H(c)=\emptyset$. Consequently, $2\Delta-2\le d_H(b)+d_H(c)=|N_H(b)\cup N_H(c)|\le |V(H)|=|V(G)|-1$. Therefore, $\Delta\le \frac{|V(G)|+1}{2}$.
\end{remark}

\begin{Theorem}\label{T1.1}
For every pair of positive real numbers $\varepsilon$ and $\xi$, there exists an integer $n_0:=n_0(\varepsilon, \xi)$ such that the following holds. Let $G$ be a $\Delta$-critical graph of order $n \ge n_0$, and suppose that there exists a vertex $a \in V(G)$ such that $H:=G-a$ is triangle-free. If $\Delta(G) \ge (\frac{2}{5}+\varepsilon)n$ and $\delta(H) \ge (\frac{1}{4}+\xi)n$, then $H$ is bipartite. Consequently, $G$ is nearly bipartite and overfull.
\end{Theorem}
	
\begin{Theorem}\label{T1.2}
Let $G$ be a $\Delta$-critical graph of order $n$ and suppose that there exists a vertex $a\in V(G)$ such that $H:=G-a$ is triangle-free.  Suppose that there exists an integer $h$ satisfying $$1 \le h<\Delta\quad\text{and}\quad 2\Delta-1<h(h+1),$$ and that, with $s:=\Delta-h$, the inequalities $$4s+\delta(H)\ge 2n-1\quad\text{and}\quad s+2\delta(H) \ge n$$ hold. Then $H$ is bipartite. Consequently, $G$ is nearly
bipartite and overfull.
\end{Theorem}
	
\begin{Corollary}\label{C1.3}
Let $G$ be a $\Delta$-critical graph of order $n$ and suppose that there exists a vertex $a \in V(G)$ such that $H:=G-a$ is triangle-free.  If $\Delta \ge \frac{3n}{7}+\left\lceil \sqrt{n} \right\rceil$ and $\delta(H) \ge \frac{2n}{7}$, then $H$ is bipartite. Consequently, $G$ is nearly
bipartite and overfull.
\end{Corollary}
	
\begin{proof}
By Remark~\ref{rem:max-degree}, we have $2\Delta-1\le n$. Let $h:=\left\lceil \sqrt{n} \right\rceil$ and $s:=\Delta-h\ge \frac{3n}{7}$. Since $1\le h <\Delta$, we have $2\Delta-1\le n\le h^2 <h(h+1)$. Note that $4s+\delta(H)\ge 2n$ and $s+2\delta(H)\ge n$. Applying Theorem \ref{T1.2} completes the proof.
\end{proof}
	
\begin{Theorem}\label{T1.4}
Let $G$ be a $\Delta$-critical graph of order $n$ and suppose that there exists a vertex $a \in V(G)$ such that $H:=G-a$ is triangle-free. Let $\sigma:=\Delta-d_G(a)$. If $4\sigma+\delta(H) \ge 2n-5$ and $\sigma+2\delta(H) \ge n-1$, then $H$ is bipartite. Consequently, $G$ is nearly bipartite and overfull.
\end{Theorem}
\begin{Corollary}
Let $G$ be a $\Delta$-critical graph of order $n$ and suppose that there exists a vertex $a \in V(G)$ such that $H:=G-a$ is triangle-free. If $\Delta \ge \frac{3n}{7}+\min\{d_G(a), \left\lceil \sqrt{n} \right\rceil\}$ and $\delta(H) \ge \frac{2n}{7}$, then $H$ is bipartite. Consequently, $G$ is nearly bipartite and overfull.
\end{Corollary}
\begin{proof}
If $\min\{d_G(a), \left\lceil \sqrt{n} \right\rceil\}=\left\lceil \sqrt{n} \right\rceil$, then applying Corollary \ref{C1.3} completes the proof. If $\min\{d_G(a), \left\lceil \sqrt{n} \right\rceil\}=d_G(a)$, then let $\sigma:=\Delta-d_G(a)\ge \frac{3n}{7}$. We have $4\sigma+\delta(H)\ge 2n$ and $\sigma+2\delta(H)\ge n$. Applying Theorem  \ref{T1.4} completes the proof.
\end{proof}
	
\section{Definitions and preliminary results}

Let $G$ be a graph with an edge $e \in E(G)$, and let $\varphi \in \mathcal{C}^k(G-e)$ be an edge coloring for some integer $k \ge \Delta$. For $u \in V(G)$, let $\varphi(u)$ be the set of colors assigned to the edges incident with $u$, and let $\overline{\varphi}(u)=[k] \setminus \varphi(u)$ be the set of colors absent from the edges incident with $u$. We call $\varphi(u)$ the set of colors {\bf present at} $u$ and $\overline{\varphi}(u)$ the set of colors {\bf missing at} $u$. Consequently, $|\varphi(u)|+|\overline{\varphi}(u)|=k$ for every vertex $u \in V(G)$. For a subset $X \subseteq V(G)$, define $\overline{\varphi}(X)=\bigcup_{u \in X} \overline{\varphi}(u)$. A set $X$ is {\bf elementary}  with respect to $\varphi$ (or $\varphi$-elementary) if $\overline{\varphi}(u) \cap \overline{\varphi}(v)=\emptyset$ for any two distinct vertices $u, v \in X$.

\begin{Definition}[Multi-fan]
Let $G$ be a graph with an edge $e=xy_1$, and an edge coloring $\phiv\in \CC^{k}(G-e)$ for some integer $k\ge\Delta$. A multi-fan at $x$ with respect to $e$ and $\phiv$ is a sequence $F=(x,e_1,y_1,\dots,e_p,y_p)$ with $p \ge 1$ consisting of edges $e_1, \dots,e_p$ and vertices $x,y_1, \dots ,y_p$ satisfying the following two conditions:
\begin{itemize}
\item [{F1.}] The edges $e_1, \dots, e_p$ are distinct, $e_1=e$, and $e_i =xy_i$ for $1\le i\le p$.
\item [{F2.}] For every edge $e_i$ with $2 \le i \le p$, there exists a vertex $y_j$ with $1 \le j \le i-1$ such that $\phiv(e_i) \in \overline{\varphi}(y_j)$.
\end{itemize}
\end{Definition}

\begin{Definition}[Kierstead path]
Let $G$ be a graph with an edge $e=y_0y_1$, and an edge coloring $\phiv\in \CC^{k}(G-e)$ for some integer $k\ge\Delta$.	
A Kierstead path with respect to $e$ and $\phiv$ is a sequence $K=(y_0,e_1,y_1,\dots,e_p,y_p)$ with $p \ge 1$ consisting of edges $e_1, \dots,e_p$ and vertices $y_0,y_1, \dots ,y_p$ satisfying the following two conditions:
\begin{itemize}
\item [{K1.}] The vertices $y_0, \dots, y_p$ are distinct, $e_1=e$, and $e_i=y_{i-1}y_i$ for $1 \le i \le p$.
\item [{K2.}] For every edge $e_i$ with $2 \le i \le p$, there exists a vertex $y_j$  with $0 \le j \le i-2$ such that $\phiv(e_i) \in \overline{\varphi}(y_j)$.
\end{itemize}
\end{Definition}

\begin{Lemma}[{\cite[Theorem 2.1]{SSTF2012}}]\label{L2.1}
Let $G$ be a graph of Class $2$ with a critical edge $e_1=xy_1$. Let $\phiv\in \CC^{\D}(G-e_1)$ and $F=(x,e_1,y_1, \dots, e_p,y_p)$ be a multi-fan at $x$ with respect to $e_1$ and $\phiv$. Then $V(F)$ is $\phiv$-elementary.
\end{Lemma}

\begin{Lemma}[{\cite[Theorem 3.3]{SSTF2012}}]\label{L2.2}
Let $G$ be a graph of Class $2$ with a critical edge $e_1=y_0y_1$. Let $\phiv\in \CC^\Delta(G-e_1)$ and $K=(y_0,e_1,y_1,e_2,y_2,e_3,y_3)$ be a Kierstead path with respect to $e_1$ and $\phiv$. Then the following statements hold.
\begin{itemize}
\item [\textup{(1)}] If $\min\{d_G(y_1),d_G(y_2)\}<\Delta(G)$, then $V(K)$ is $\phiv$-elementary.
\item [\textup{(2)}] $|\overline{\varphi}(y_3) \cap (\overline{\varphi}(y_0) \cup \overline{\varphi}(y_1))| \le 1$.	\end{itemize}
\end{Lemma}

Throughout this paper, all multigraphs considered are
loopless, although parallel edges are allowed. A multigraph is called nearly bipartite if deleting at most one vertex results in a bipartite multigraph.

\begin{Lemma}[{\cite[Theorem 2]{LM1986}}]\label{L2.3}
Let $G$ be a nearly bipartite multigraph. Then $\chi'(G)=\Delta(G)$ if and only if $G$ contains no $\Delta(G)$-overfull subgraph.
\end{Lemma}

\begin{Lemma}\label{L2.8}
Let $G$ be a  $\Delta$-critical graph. If there exists a vertex $v \in V(G)$ such that $G-v$ is bipartite, then $G$ is overfull.
\end{Lemma}
\begin{proof}
Since $G$ is $\Delta$-critical, $\chi'(G)=\Delta+1$. Moreover, $G-v$ is bipartite, so $G$ is nearly bipartite. By Lemma \ref{L2.3}, $G$ contains a $\Delta$-overfull subgraph $G_1$. If $G_1$ is a proper subgraph of $G$, then $\chi'(G_1) \le \Delta$. On the other hand, $|E(G_1)|>\Delta\lfloor \frac{|V(G_1)|}{2} \rfloor$ implies that $G_1$ cannot be properly edge-colored with $\Delta$ colors. Hence $\chi'(G_1)>\Delta$, a contradiction. Therefore, $G_1=G$ and $G$ is overfull. This completes the proof.
\end{proof}

\begin{Lemma}[{\cite{LARA2020}}]\label{L2.4}
Let $G$ be a $\Delta$-critical graph. For every vertex $u\in V(G)$, there exist  at least two distinct vertices $v_1, v_2\in N_G(u)$ such that $d_G(v_i)=\Delta$ and $\sum_{w \in N_G(v_i)} d_G(w) \ge \Delta^2 -\Delta +1$ for every $i\in \{1,2\}$. Equivalently, $\sum_{w \in N_G(v_i)}(\Delta -d_G(w)) \leq \Delta-1$ for every $i\in \{1,2\}$.
\end{Lemma}
	
\begin{Lemma}[{\cite{CC2020}}]\label{L2.5}
For every $\varepsilon \in (0,1)$, there exists a positive constant $D(\varepsilon)$ such that every $\Delta$-critical graph $G$ with $\Delta\ge D(\varepsilon)$ satisfies $\overline{d}(G) \ge (1-\varepsilon)\Delta$.
\end{Lemma}
The following lemma corresponds to the $r=3$ case of Theorem~1.1 in Andr{\'a}sfai, Erd{\H{o}}s, and S{\'o}s~\cite{AES1974}, originally established as their Lemma~1.2.
\begin{Lemma}[\cite{AES1974}]\label{L2.6}
If $G$ is a triangle-free graph of order $n$ with minimum degree $\delta(G)>\frac{2n}{5}$, then $G$ is bipartite.
\end{Lemma}
	
\begin{Lemma}[Bipartite extension lemma]\label{L2.9}
Let $G$ be a triangle-free graph, and let $X$ and $Y$ be two disjoint independent sets. Let $Z:=V(G) \setminus (X\dot{\cup} Y)$, $x:=|X|$, $y:=|Y|$,  and $z:=|Z|$, where $\dot{\cup}$ denotes the disjoint union.
Suppose that real numbers $d_1, d_2$ satisfy
$d_G(u) \ge d_1$ for every $u \in X$ and  $d_G(v) \ge d_2$ for every $v \in Y$.
If $$\delta(G)>z-1, \,\,  \delta(G)+d_1+d_2>x+y+2z, \,\, and \,\, 2\delta(G)>\max\{x, y\}+z,$$
then $G$ is bipartite.
\end{Lemma}

\begin{proof}
If $Z=\emptyset$, then the assertion follows from the independence of $X$ and $Y$. Assume that $Z \neq \emptyset$. We claim that for every $w\in Z$, exactly one of $N_G(w)\cap X$ and $N_G(w)\cap Y$ is nonempty. We rule out the following two possibilities. First, if $N_G(w)\cap (X\dot{\cup} Y)=\emptyset$, then $d_G(w) \le z-1<\delta(G)$, contradicting the fact that $d_G(w)\ge \delta(G)$. Next, suppose that $w$ has a neighbor $u \in X$ and a neighbor $v \in Y$. Let $A:=|N_G(w) \cap X|$, $B:=|N_G(w) \cap Y|$, and $C:=|N_G(w) \cap Z|$. Since $X$ is independent, $N_G(u) \subseteq Y\dot{\cup} Z$. Moreover, since $G$ is triangle-free and $uw\in E(G)$, the sets $N_G(u)$  and  $N_G(w)\cap (Y\dot{\cup} Z)$ are disjoint. Hence, we get $d_G(u)\le y+z-B-C$. Similarly, $d_G(v)\le x+z-A-C$. Since $d_G(w)=A+B+C$, we have $d_G(u)+d_G(v)+d_G(w)\le x+y+2z-C$. It follows that $d_G(w)\le x+y+2z-d_1-d_2-C<\delta(G)$, contradicting the fact that $d_G(w)\ge \delta(G)$.

Next, define $Z_X:=\{w \in Z:N_G(w) \cap X=\emptyset\}$ and $Z_Y:=\{w \in Z:N_G(w)\cap Y=\emptyset\}$.
Observe that $Z=Z_X \dot{\cup} Z_Y$. Furthermore, there are no edges within $X$ or $Y$, nor between $X$ and $Z_X$, nor between $Y$ and $Z_Y$. We claim that both $Z_X$ and $Z_Y$ are  independent sets. Otherwise, without loss of generality, suppose that $Z_X$ is not an independent set.
There exist adjacent vertices $w_1$ and  $w_2$ in $Z_X$. Note that $N_G(w_1), N_G(w_2) \subseteq Y \dot{\cup} Z$. Since $G$ is triangle-free, $N_G(w_1)\cap N_G(w_2)=\emptyset$. It follows that $2\delta(G)\le d_G(w_1)+d_G(w_2)=|N_G(w_1)\cup N_G(w_2)|\le y+z<2\delta(G)$, which is a contradiction. The same argument shows that $Z_Y$ is independent, using $x+z<2\delta(G)$. Consequently, $X\dot{\cup} Z_X$ and $Y\dot{\cup} Z_Y$ constitute a bipartition of $G$, which proves that $G$ is bipartite. This completes the proof.
\end{proof}

\section{Proof of Theorem 1.1}
\begin{proof}
Let $p:=\min\{\frac{\varepsilon}{4}, \frac{\xi}{4}, \frac{1}{100}\}$ and $q:=\frac{p^2}{2}$. By Lemma \ref{L2.5}, there exists a positive constant $D(q)$ such that every $\Delta$-critical graph with maximum degree at least $D(q)$ has average degree at least $(1-q)\Delta$. Choose an integer $n_0=n_0(\varepsilon,\xi)$ sufficiently large such that $n_0\ge \max\{\frac{2}{p^2}, \frac{5}{2}D(q),3\}$. Since $n\ge n_0$, we have $\Delta \ge (\frac{2}{5}+\varepsilon)n>\frac{2}{5}n_0\ge D(q)$. Given that $G$ is $\Delta$-critical, Lemma \ref{L2.5} yields $\overline{d}(G) \ge (1-q)\Delta$.
		
Define $D_G:=\sum_{u \in V(G)}(\Delta -d_G(u))=n(\Delta-\overline{d}(G)) \le nq\Delta \le qn^2$. Let $D_H:=\sum_{u \in V(H)}(\Delta-d_H(u))$. Relating the two quantities yields $D_H=D_G-(\Delta-d_G(a))+d_G(a)=D_G-\Delta+2d_G(a) \le D_G+\Delta \le qn^2+n$. Since $n \ge \frac{2}{p^2}$, we have $n \le \frac{p^2n^2}{2}$. Consequently, $D_H \le qn^2+n \le\frac{p^2}{2}n^2+\frac{p^2n^2}{2}=p^2 n^2$.
Now, define $Z:= \{u \in V(H) : d_H(u)<\Delta-pn\}$ and $z:=|Z|$. We claim that $z<pn$. If $Z=\emptyset$, then the inequality holds trivially. If $Z\neq \emptyset$, then we have $zpn<\sum_{u \in Z} (\Delta-d_H(u))\le D_H\le p^2 n^2$. Therefore, we get $z<pn$.

Let $H':=H-Z$. Clearly, $V(H')\neq \emptyset$ since $|V(H')|=n-1-z>(1-p)n-1>0$. For every $u \in V(H')$, we have $d_H(u) \ge \Delta-pn$. Since deleting the vertices in $Z$ removes at most $z$ neighbors from any remaining vertex, we get $d_{H'}(u) \ge \Delta-pn-z>\Delta-2pn$. Combining $\Delta \ge (\frac{2}{5}+\varepsilon)n$ and $p \le \frac{\varepsilon}{4}$ yields $d_{H'}(u)>\Delta -2p n \ge (\frac{2}{5}+\frac{\varepsilon}{2})n>\frac{2}{5}n \ge \frac{2}{5}|V(H')|$. Therefore, $\delta(H')>\frac{2}{5}|V(H')|$. Since $H$ is triangle-free and $H' \subseteq H$, $H'$ is also triangle-free. By Lemma \ref{L2.6}, $H'$ is bipartite. Let $V(H')=X\dot{\cup} Y$, where $X$ and $Y$ are the two partite sets of $H'$. Since $\delta(H')>0$, $X$ and $Y$ are nonempty. Let $r:=\Delta-2pn$. For every $u \in V(H')$, $r<d_{H'}(u) \le \Delta$. Note that $|E(H')|=\sum_{u \in X} d_{H'}(u)=\sum_{u \in Y} d_{H'}(u)$. Counting from $X$, $|E(H')|>r|X|$. Counting from $Y$, $|E(H')|\le \Delta|Y|$. Using $|Y|=|V(H')|-|X|$, it follows that $r|X|<\Delta(|V(H')|-|X|)$, which gives $|X|<\frac{\Delta}{\Delta+r}|V(H')|$. Substituting $\Delta >\frac{2n}{5}$ and $p\le \frac{1}{100}$ into the above inequality gives $|X|<\frac{\Delta}{\Delta+r}|V(H')| \le \frac{\Delta n}{2\Delta-2pn}< \frac{n}{2-5p} \le (\frac{1}{2}+2p)n$. Similarly, $|Y|<(\frac{1}{2}+2p)n$.  Consequently, there exist two disjoint independent sets $X$ and $Y$ such that $|X|, |Y|<(\frac{1}{2}+2p)n$ and $d_H(u)\ge \Delta-pn$ for every $u\in X \dot{\cup} Y$. Since $V(H')=X\dot\cup Y$ and $H'=H-Z$,
we have $V(H)=X\dot\cup Y\dot\cup Z$, and $|Z|=z<pn$.

We now apply Lemma~\ref{L2.9} to $H$ with the independent sets $X$ and $Y$. By the construction of $X$ and $Y$, we have  $x,y<(\frac{1}{2}+2p)n$ and may take $d_1=d_2=\Delta-pn$. Note that $z=n-1-(x+y)$. Substituting these into the conditions of Lemma \ref{L2.9}, it suffices to verify that $\delta(H)>z-1$, $n-1+z-2\Delta+2pn<\delta(H)$, and $2\delta(H)>\max\{x,y\}+z$. The first inequality follows from $\delta(H)\ge(\frac{1}{4}+\xi)n>pn>z$. For the second inequality, using $z<pn$ and $\Delta \ge (\frac{2}{5}+\varepsilon)n$, we obtain $n-1+z-2\Delta+2pn<(\frac{1}{5}-2\varepsilon+3p)n-1\le(\frac{1}{5}-\frac{5}{4}\varepsilon)n-1<\delta(H)$, where we have used $p\le \frac{\varepsilon}{4}$. Finally, since $x,y<(\frac{1}{2}+2p)n$ and $p\le \frac{\xi}{4}$, we have $\max\{x,y\}+z<(\frac{1}{2}+3p)n\le(\frac{1}{2}+\frac{3\xi}{4})n<2\delta(H)$. Hence $H$ is bipartite, and Lemma~\ref{L2.8} implies that $G$ is overfull. This completes the proof.
\end{proof}

\section{Proof of Theorem 1.2}
\begin{proof}
Applying Lemma \ref{L2.4} to $a$,  choose a vertex $b\in N_G(a)$ satisfying $d_G(b)=\Delta$ and $\sum_{u \in N_G(b)} (\Delta-d_G(u)) \le \Delta-1$. Applying Lemma~\ref{L2.4} to $b$, we obtain at least two distinct neighbors satisfying the stated properties. Since at most one of them is $a$, we may choose $c\in N_G(b)\setminus\{a\}$ such that $d_G(c)=\Delta$
and $\sum_{u \in N_G(c)} (\Delta-d_G(u)) \le \Delta-1$. Note that $d_H(b)=\Delta-1$ and $d_H(c)\ge\Delta-1$.

We claim that there exist two disjoint independent sets $X$ and $Y$ satisfying $X \subseteq N_H(b)$, $Y \subseteq N_H(c)$, $|X|=|Y|=s$, and $d_H(u)\ge s$ for every $u \in X\dot{\cup} Y$. It suffices to consider $N_H(b)$, as the argument for $N_H(c)$ is identical. Let $L_b:=\{u \in N_H(b):d_H(u) \le \Delta-h-1\}$. For every $u \in L_b$, $\Delta-d_H(u) \ge h+1$. Therefore, $\sum_{u \in L_b} (\Delta-d_H(u)) \ge (h+1)|L_b|$. Since $L_b \subseteq N_H(b)$, we have $\sum_{u \in L_b} (\Delta-d_H(u)) \le \sum_{u\in N_H(b)} (\Delta-d_H(u)) =\sum_{u\in N_H(b)} (\Delta-d_G(u))+|N_H(b)\cap N_G(a)|\le(\Delta-1)+d_G(a)\le 2\Delta-1$. Thus, $(h+1)|L_b| \le 2\Delta-1$, which implies $|L_b| \le \lfloor \frac{2\Delta-1}{h+1} \rfloor <h$. Therefore, $|N_H(b) \setminus L_b|\ge (\Delta-1)-(h-1)=\Delta-h=s$. We may therefore choose $s$ vertices from $N_H(b) \setminus L_b$ to form $X$. Since $H$ is triangle-free, $X$ is an independent set. Now we get an independent set $X \subseteq N_H(b)$ with $|X|=s$ and $d_H(u) \ge s$ for every $u \in X$. Applying the same argument to $N_H(c)$, we obtain an independent set $Y \subseteq N_H(c)$ with $|Y|=s$ and $d_H(u)\ge s$ for every $u\in Y$.  Since $H$ is triangle-free and $bc\in E(H)$, it follows that  $X \cap Y=\emptyset$.

Let $Z:=V(H)\setminus (X\dot{\cup} Y)$ and $z:=|Z|=n-1-2s$. We apply Lemma~\ref{L2.9} to $H$ with the independent sets $X$ and $Y$. By the construction of $X$ and $Y$, we may take $x=y=s$ and $d_1=d_2=s$. Thus, the three conditions of Lemma \ref{L2.9} reduce to $\delta(H)>z-1$, $\delta(H)>2z$, and $2\delta(H)>s+z$. Since $z\ge 0$, the second inequality implies the first inequality. Hence, it suffices to verify the latter two inequalities. By the hypothesis of the theorem and the identity $z=n-1-2s$, we have $\delta(H)\geq2n-4s-1=2z+1>2z$. Moreover, $2\delta(H)\geq n-s=s+z+1>s+z$. Hence $H$ is bipartite, and Lemma~\ref{L2.8} implies that $G$ is overfull. This completes the proof.
\end{proof}

\section{Proof of Theorem 1.4}
\begin{proof}
When $\Delta=2$, the $\Delta$-critical graph $G$
is connected and satisfies $\chi'(G)=3$.
A connected graph of maximum degree $2$ is either a
path or a cycle, and the condition $\chi'(G)=3$
forces $G$ to be an odd cycle. Deleting any vertex results in a path that is bipartite. Moreover, $G$ is overfull since $|E(G)|=n>n-1=2 \cdot \lfloor \frac{n}{2} \rfloor$. Now assume that $\Delta \ge 3$. By Remark~\ref{rem:max-degree}, there exists a vertex $b \in N_G(a)$ such that $d_G(b)=\Delta$. Thus, $d_H(b)=\Delta -1$ and $2\Delta-1\le n$.
We claim that $\Delta \ge 2d_G(a)-1$. Otherwise, we have $\Delta\le 2d_G(a)-2$. Since $4\sigma+\delta(H) \ge 2n-5$ and $2\Delta-1\le n$, it follows that $\delta(H) \ge 2n-4\sigma-5 \ge 4d_G(a)-7 \ge 2\Delta-3>\Delta-1$, which contradicts $\delta(H) \le d_H(b)=\Delta-1$.
	
Let $\phiv \in \mathcal{C}^\Delta(G-ab)$. Since $(a,ab,b)$ is a multi-fan with respect to $ab$ and $\phiv$, $\{a, b\}$ is $\phiv$-elementary by Lemma \ref{L2.1}. Then we have $|\overline{\phiv}(a) \cup \overline{\phiv}(b)|=\Delta-d_G(a)+2$ and $|[\Delta] \setminus (\overline{\phiv}(a) \cup \overline{\phiv}(b))|=d_G(a)-2$. Let $X:=\{u \in N_G(b) \setminus \{a\}:\phiv(bu) \in \overline{\phiv}(a)\}$. Since $\{a, b\}$ is elementary, $|X|=|\overline{\phiv}(a)|=\Delta-d_G(a)+1$. Order the vertices of $X$ arbitrarily as $u_1,\ldots,u_{|X|}$. Since $\varphi(bu_i)\in\overline{\varphi}(a)$ for every $i\in[|X|]$, the sequence $(b,ba,a,bu_1,u_1,\ldots,bu_{|X|},u_{|X|})$ is a multi-fan with respect to $ba$ and $\varphi$. By Lemma \ref{L2.1}, $\{a, b\}\cup X$ is $\phiv$-elementary. This yields $\sum_{u \in X} (\Delta-d_G(u))\le \Delta-|\overline{\varphi}(a)\cup\overline{\varphi}(b)|=d_G(a)-2$. Consequently, for every $u\in X$, we have $d_G(u)\ge \Delta-d_G(a)+2$ and $d_H(u)\ge \Delta-d_G(a)+1$. Since $\Delta\ge 2d_G(a)-1$, we obtain $|X|=\Delta-d_G(a)+1 \ge d_G(a)$. Given that $|N_G(a) \setminus \{b\}|=d_G(a)-1$, we can find a vertex $u_0 \in X$ such that $au_0 \notin E(G)$.
	
Set $Y:=\{v \in N_G(u_0):\phiv(u_0v) \in \overline{\phiv}(a) \cup \overline{\phiv}(b)\}$. Since $\overline{\phiv}(u_0) \cap (\overline{\phiv}(a) \cup \overline{\phiv}(b))=\emptyset$, every color in $\overline{\phiv}(a) \cup \overline{\phiv}(b)$ is present at $u_0$. Consequently, $|Y|=\Delta-d_G(a)+2$. For every $v \in Y\setminus \{b\}$, the four vertices $a,b,u_0,v$ are pairwise distinct.
Since $\phiv(bu_0)\in \overline{\phiv}(a)$ and $\phiv(u_0v)\in \overline{\phiv}(a) \cup \overline{\phiv}(b)$, the sequence $(a,ab,b,bu_0,u_0,u_0v,v)$ forms a Kierstead path. By Lemma \ref{L2.2}(2), $|\overline{\phiv}(v) \cap (\overline{\phiv}(a) \cup \overline{\phiv}(b))| \le 1$, which implies $|\overline{\phiv}(v)|\le d_G(a)-1$, and hence $d_G(v)\ge \Delta-d_G(a)+1$. Thus, $d_H(v) \ge \Delta-d_G(a)$. Since $u_0 \in X$, we have $\varphi(u_0b) \in\overline{\phiv}(a)$, and therefore $b \in Y$. Noting that $d_H(b)=\Delta-1\ge\Delta-d_G(a)$, it follows that $d_H(v)\ge \Delta-d_G(a)$ for every $v\in Y$. Since $H$ is triangle-free and $bu_0 \in E(H)$, the sets $N_H(b)$ and $N_H(u_0)$ are independent and disjoint. Hence $X$ and $Y$ are disjoint independent sets. Note that $\sigma=\Delta-d_G(a)$. Consequently, there exist two disjoint independent sets $X$ and $Y$ such that $|X|=\sigma+1$, $|Y|=\sigma+2$, $d_H(u) \ge \sigma+1$ for every $u \in X$, and $d_H(v)\ge \sigma$ for every $v \in Y$.
	
Let $Z:=V(H)\setminus (X\dot{\cup} Y)$ and $z:=|Z|=n-2\sigma-4$.  We now apply Lemma~\ref{L2.9} to $H$ with the independent sets $X$ and $Y$. By the construction of $X$ and $Y$, we may take $x=\sigma+1$, $y=\sigma+2$, $d_1=\sigma+1$, and $d_2=\sigma$. Thus, the three conditions of Lemma \ref{L2.9} reduce to $\delta(H)>z-1$, $\delta(H)>2z+2$, and $2\delta(H)>\sigma+z+2$. Since $z\ge 0$, the second inequality implies the first inequality. Hence, it suffices to verify the latter two inequalities. By the hypothesis of the theorem and the identity $z=n-2\sigma-4$, we have $\delta(H)\geq2n-4\sigma-5=2z+3$ and  $2\delta(H)\geq n-\sigma-1=\sigma+z+3$. Hence $H$ is bipartite, and Lemma~\ref{L2.8} implies that $G$ is overfull. This completes the proof.
\end{proof}

\noindent {\bf Acknowledgements}

This work was supported by Hebei Natural Science Foundation (A2026208005) and Foundation of China Scholarship Council (202508130088).

\end{document}